\documentclass[10pt,a4paper]{article}
\usepackage{amsmath,latexsym,amsfonts,amssymb}
\usepackage{wrapfig}
\usepackage{amscd,array,epic,eepic,calc,bbm,float}
\usepackage{ifthen,psfrag,verbatim,epsfig,graphicx,enumerate,subfig}
\usepackage{makeidx}
\usepackage{amsthm}
\usepackage{xcolor}

\def\bql{\begin{equation}\label}
\def\eql{\end{equation}\noindent}

\def\brl{\begin{eqnarray}\label}
\def\erl{\end{eqnarray}\noindent}

\def\bro{\begin{eqnarray*}}
\def\ero{\end{eqnarray*}\noindent}

\def\brr{\begin{array}}
\def\err{\end{array}\noindent}

\def\bdl{\begin{display}\label}
\def\edl{\end{display}\noindent}

\def\bdo{\begin{display}}
\def\edo{\end{display}\noindent}

\def\bth{\begin{theorem}}
\def\eth{\end{theorem}}

\def\bcr{\begin{corollary}}
\def\ecr{\end{corollary}}

\def\bpr{\begin{proposition}}
\def\epr{\end{proposition}}

\def\blm{\begin{lemma}}
\def\elm{\end{lemma}}

\def\bdf{\begin{definition}}
\def\edf{\end{definition}}

\def\bas{\begin{assumptions}}
\def\eas{\end{assumptions}}

\def\bex{\begin{example}\rm}
\def\eex{\end{example}}

\def\bxx{\begin{exercise}\rm}
\def\exx{\end{exercise}}

\def\brm{\begin{remark}\rm}
\def\erm{\end{remark}}

\def\bma{\begin{pmatrix}}
\def\ema{\end{pmatrix}}

\def\bcs{\begin{cases}}
\def\ecs{\end{cases}}

\def\btb{\begin{center}\begin{tabular}}
\def\etb{\end{tabular}\end{center}}

\def\bit{\begin{itemize}}
\def\eit{\end{itemize}}

\def\df{\par\noindent{\bf Definition }}

\def\bew{\par\noindent{\bf Proof }}
\def\qed{\quad\hfill\mbox{\P}}
\def\qend{\hfill\mbox{$\lozenge$}}

\def\c{\chi}
\def\d{\delta}
\def\e{\epsilon}

\def\f{\varphi}
\def\g{\gamma}
\def\h{\eta}

\def\j{\psi}

\def\l{\lambda}
\def\m{\mu}
\def\n{\nu}
\def\p{\pi}
\def\q{\theta}
\def\r{\rho}

\def\u{\upsilon}
\def\w{\omega}
\def\x{\xi}

\def\AC{{\cal A}}

\def\DC{{\cal D}}

\def\HC{{\cal H}}

\def\LB{{\bf L}}

\def\XC{{\cal X}}

\def\UB{{\mathbf U}}

\def\WB{{\mathbf W}}
\def\XB{{\mathbf X}}
\def\YB{{\mathbf Y}}
\def\ZB{{\mathbf Z}}

\def\ab{{\mathbf a}}
\def\bb{{\mathbf b}}
\def\eb{{\mathbf e}}

\def\pb{{\mathbf p}}
\def\qb{{\mathbf q}}
\def\ub{{\mathbf u}}

\def\wb{{\mathbf w}}
\def\xb{{\mathbf x}}

\def\zb{{\mathbf z}}
\def\oneb{{\bf 1}}
\def\zerob{{\bf0}}
\def\nfb{\boldsymbol{\nf}}

\def\alb{\boldsymbol{\alpha}}
\def\betab{\boldsymbol{\beta}}

\def\lp{\ell_p}

\def\pbb{{\mathbb P}}
\def\rbb{{\mathbb R}}

\def\imp{\Rightarrow}

\def\inv{^{-1}}
\def\ginv{^{\leftarrow}}

\def\nf{\infty}

\def\ov{\overline}
\def\prl{\partial}
\def\sm{\setminus}
\def\ss{\subset}

\def\nt{\noindent}

\newtheorem{theorem}{Theorem}[section]
\newtheorem{lemma}[theorem]{Lemma}
\newtheorem{proposition}[theorem]{Proposition}
\newtheorem{corollary}[theorem]{Corollary}
\newtheorem{definition}{Definition}
\newtheorem{example}{Example}
\newtheorem{exercise}{Exercise}
\newtheorem{remark}{Remark}
\newtheorem{assumptions}{Assumptions}[section]

\makeindex

\begin{document}
\date{\today}
\title{Asymptotic independence for unimodal densities}
\author{\renewcommand{\thefootnote}{\arabic{footnote}}Guus Balkema\footnotemark[1] $\qquad$ Natalia Nolde\footnotemark[2]}
\footnotetext[1]{Department of Mathematics, University of Amsterdam, Science Park 904, 1098XH Amsterdam, Netherlands\\ {\it Email address:} A.A.Balkema@uva.nl} \footnotetext[2]{Department of Mathematics and RiskLab, ETH Zurich, Raemistrasse 101, 8092 Zurich, Switzerland\\ {\it Email address:} natalia.lysenko@math.ethz.ch} \addtocounter{footnote}{2}%
\maketitle
\setcounter{page}{1}

\begin{abstract}
Asymptotic independence of the components of random vectors is a concept used in  many applications. 
The standard criteria for checking asymptotic independence are given in terms of distribution functions (dfs). Dfs are rarely available in an explicit form, especially in the multivariate case. Often we are given the form of the density or, via the shape of the data clouds, one can obtain a good geometric image of the asymptotic shape of the level sets of the density. This paper establishes a simple sufficient condition for asymptotic independence for light-tailed densities in terms of this asymptotic shape. This condition extends Sibuya's classic result on asymptotic independence for Gaussian densities. 
\end{abstract}

\nt{\bf Key words:} Asymptotic independence; blunt; homothetic density; level set; skew-normal; star-shaped
\nt{\bf 2000 MSC:} 60G55, 60G70, 62E20

\setcounter{equation}{0}
\setcounter{section}{0}

\section{Introduction}

The purpose of the present paper is to provide simple sufficient conditions that  ensure asymptotic independence of the components of random vectors whose probability distribution is described by a density. Standard criteria for checking asymptotic independence are given in terms of distribution functions (dfs). However, these are not always available in an explicit form in the multivariate case, and they give little insight in what large samples from a distribution will look like. Often we are given  a density in analytic form. For light-tailed densities, the data clouds give a good geometric image of the asymptotic shape of the level sets. Hence it is of interest to have conditions for asymptotic independence in terms of the shape of the level sets of the underlying density, or in terms of a limiting shape for data clouds. 

For vector valued data it is standard practice to plot the bivariate sample clouds for all component pairs. In our final result it is the asymptotic behaviour of the shape of these bivariate sample clouds, as the size of the data set increases, that determines asymptotic independence of the coordinates for the underlying multivariate distribution. 

Unimodal densities whose level sets all have the same shape are called homothetic. The decay along any ray then is the same up to a scale constant depending on the direction, and hence the concept of light and heavy tails is well-defined. This remains true if we only assume that the level sets have the same shape asymptotically. Our primary focus here is on light-tailed densities, but for additional insight we also include some results for the heavy-tailed case. Throughout the paper, we assume continuity of dfs and of densities. 

Intuitively, for bivariate data, asymptotic independence means that large values in one coordinate are unlikely to be accompanied by large values in the other coordinate. Situations with a low chance of simultaneous extremes are often encountered in practice, for example in applications which involve modeling environmental data (e.g. \cite{Ledford1996}) or network traffic data (e.g. \cite{Maulik2002}); see \cite{Resnick2002} for further references. It is a well-known result, dating back to 1959 (see \cite{Sibuya1960}), that the components of a vector with a Gaussian density are asymptotically independent whatever the correlation. Asymptotic independence also holds for light-tailed elliptical densities (see e.g. \cite{Hult2002}). As a further generalization we shall show that a vector with a continuously differentiable homothetic light-tailed density whose level sets are convex has asymptotically independent components. 

\bth\label{thmai0} Let $D$ be a bounded open convex set in $\rbb^d$, containing the origin, with a $C^1$ boundary (i.e., at each boundary point there is a unique tangent plane). Let $c_n>0$ decrease to zero such that $c_{n+1}/c_n\to0$ and let $r_n$ be positive reals such that $r_{n+1}/r_n\to1$. If $\XB$ is a random vector in $\rbb^d$ with a continuous probability density $f$ whose level sets satisfy \bql{q1t1}\{f>c_n\}=r_nD\qquad n\ge n_0,\eql
then the components of $\XB$ are asymptotically independent. 
\eth
We shall prove more general results, Theorem~\ref{thmai3} and \ref{thmai1} below. The generalizations we introduce are simple, but they result in theorems for which some extra terminology has to be developed. 
\bit
\item The condition that the level sets all be of a given shape is replaced by the condition that the level sets can be scaled to converge to a limit shape $D$. Theorem~\ref{thmai0} remains valid if we replace \eqref{q1t1} by the limit relation $\{f>c_n\}/r_n\to D$.
\item Asymptotic independence holds if it holds for the bivariate marginals; in this spirit, our theorem imposes conditions on the projection $D_{12}$ of the set $D$ on the $x_1,x_2$-coordinate plane.
\item Convex level sets are replaced by star-shaped level sets.
\item The condition of a smooth ($C^1$) boundary is replaced by a condition which only affects the maximum: the coordinatewise supremum of the points in $D_{12}$ should not be a boundary point of $D_{12}$. Such a set is called blunt.
\item Our final result (Theorem~\ref{thmai2}) is on bivariate sample clouds. If these converge onto a blunt star-shaped set, asymptotic independence holds.
\eit   
The conditions in Theorem~\ref{thmai0} above on the sequences $(c_n)$ and $(r_n)$ ensure that the density has light tails. If $c_{n+1}/c_n\to1$ then the density has heavy tails, and the same condition on the shape of the level sets implies asymptotic dependence. 

The paper is organized as follows. In Section~2 we discuss the concept of asymptotic independence. Section~3 introduces the class $\HC$ of continuous multivariate homothetic densities with star-shaped level sets. The simple structure of these densities makes them a good starting point in our investigation of the relation between the shape of the level sets and asymptotic independence. Here we formulate our main results (Theorem~\ref{thmai3}~and~\ref{thmai1}) for light-tailed densities, and a counterpart (Theorem~\ref{thmad1}) for heavy-tailed densities. In Section~4 we introduce sample clouds, and state conditions for asymptotic independence in terms of their limiting shape (Theorem~\ref{thmai2}). Section~5 provides various examples. The appendix clears up a number of minor points in the main text, by providing supplementary results and counterexamples.

\section{Asymptotic independence}

In this section we discuss asymptotic independence for multivariate distributions. We begin with a heuristic approach in the bivariate setting. This will clarify the significance of the concept for risk management. We discuss the commonly used criteria for asymptotic independence, the relation with multivariate extreme value theory, and describe possible forms of asymptotic dependence. 

\subsection{Heuristics}

In finance one is interested in the future value of stocks, say in one year's time. Let $X$ and $Y$ denote the future value of two stocks. Suppose the distribution of the pair $(X,Y)$ is given by a df $F$ with continuous marginals $F_1$ and $F_2$. One is concerned about the risk that the stocks have a low value at this future date. Let $x_p$ denote the $p$-quantile of the first stock, $X$, and $y_p$ the $p$-quantile of the second. The probability that the values of both stocks lie below the $p$-quantile for some given small value of $p\in(0,1)$ is $F(x_p,y_p)$.  A risk averse investor would like this probability (of simultaneous loss) to be small compared to the probability $p$ of a loss in either of the stocks. Future extreme low values of the two stocks are said to be asymptotically independent if $F(x_p,y_p)=o(p)$ for $p\to0$. Asymptotic dependence will increase risk for a portfolio containing these two stocks.

A different example of risk is presented by the yearly maxima for high water levels at particular points on the coast of Holland, Great Britain and the U.S., say IJmuiden, Harwich and New Orleans. Consider a data set of 200 observations stretching back to the beginning of the nineteenth century. We assume that the data have been standardized to offset tidal effects. Now for each of these locations pick the five largest values. This yields three subsets of five elements in the set of 200 years. One would suppose that there is considerable overlap between the five years selected for IJmuiden and for Harwich, since in both cases the cause is the same, a North Western storm in the North Atlantic forcing water into the funnel formed by the West coast of Holland and Belgium and the South coast of Great Britain, and opening into the Southern Atlantic via the Channel. High water levels in New Orleans have a different cause. So sea levels in New Orleans and in IJmuiden should be
asymptotically independent, and one would expect the corresponding five point subsets to be disjoint with high probability. In general, given a positive integer $k$ and a sequence of independent observations $\ZB_1,\ZB_2,\ldots$ from a bivariate df $F$, one could look in the sample cloud of the first $n$ points $\ZB_i=(X_i,Y_i)$ at the $k$ largest observations for the coordinate $X_i$ and for the coordinate $Y_i$. This yields two subsets of the index set $\{1,\ldots,n\}$. Let $p_n(k)$ denote the probability that the two subsets have a point in common. In the case of asymptotic independence one would expect that for fixed $k$, $p_n(k)\to0$ as $n\to\nf$. This is indeed the case, as shown in \cite{Gnedin1993} (Proposition~2).

A related way to understand extremal dependence is by looking at the probability of a \emph{record}, $p_n(1)$, the probability that the coordinatewise maximum of a sample of $n$ points is given by one of the sample points. It is shown in \cite{Gnedin1994} (Theorem~2) that this probability vanishes as $n\to\nf$ exactly when the underlying vector has asymptotically independent components. For the sake of completeness we give an alternative proof of this basic result in the Appendix; see Proposition~\ref{prec}. 


\subsection{Sibuya's condition}

In his seminal paper \cite{Sibuya1960} on multivariate extremes, Sibuya  shows that the components of bivariate vectors with normal densities, whatever their correlation, are asymptotically independent. For a vector $\ZB=(X,Y)$ with df $F$ and continuous marginals $F_1$ and $F_2$ he introduces a function $P$ by
\bql{q2PF}
P(F_1(x),F_2(y))=\pbb\{X>x, Y>y\}=1+F(x,y)-F_1(x)-F_2(y).
\eql
The function $P$ is well-defined and continuous on the unit square. We can now give \emph{Sibuya's condition} for asymptotic independence (see Theorem~2 in \cite{Sibuya1960}).

\df Let $\ZB=(X,Y)$ have df $F$ with continuous marginals. The components $X$ and $Y$ are \emph{asymptotically independent} if the function $P$ introduced above satisfies
\bql{q2ai}
P(1-s,1-s)=o(s)\qquad s>0, s\to0.
\eql
By abuse of language we also say that the vector $\ZB$ or the df $F$ is \emph{asymptotically independent}. 

Independence of a bivariate vector is not affected by the marginal distributions. Similarly for asymptotic independence. It is preserved under coordinatewise increasing transformations. Since we assume that the marginals $F_1$ and $F_2$ are continuous, there is a unique function $C$ on the unit square such that $F(x,y)=C(F_1(x),F_2(y))$. The function $C$ is known as the \emph{copula} of the df $F$. It is a df on the unit square with uniform marginals. Sibuya's condition is a condition on the copula $C$ since $P(u,v)=1+C(u,v)-u-v$ for $u,v\in[0,1]$. If the condition holds for a vector $(X,Y)$ with a continuous df $F$, it automatically holds for any vector whose df is continuous and has the same copula as $F$.

In this paper we consider asymptotic independence for maxima. For minima one would define asymptotic independence similarly, in terms of the copula, by: $C(s,s)=o(s)$ for $s\downarrow0$. 

\bpr\label{p2ab}
Suppose there exist $a,b>0$ such that
\bql{q2Pab}
P(1-as,1-bs)/s\to0\qquad s\downarrow0.
\eql
Then asymptotic independence holds, and the relation~(\ref{q2Pab}) is valid for all $a,b>0$.
\epr
\bew We may assume $a\le b$ by symmetry and $a=1$ replacing $s$ by $as$ in the denominator. By monotonicity Sibuya's condition holds. A similar argument gives~(\ref{q2Pab}) for any positive $a,b$. \qed

For a vector $(X,Y)$ with marginal dfs $F_1$ and $F_2$, Sibuya's condition may be formulated in terms of conditional quantile exceedances as
\bql{q2con}
\l_U(X,Y):=\lim_{q\uparrow1}\pbb\{X>F_1\ginv(q)\mid Y>F_2\ginv(q)\}=0,
\eql
where $F\ginv(q):=\inf\{x\in\rbb\mid F(x)\ge q\}$ for $q\in(0,1)$ denotes the (minimal) $q$-quantile of $F$. The limit $\l_U\in[0,1]$, if it exists, is known as the \emph{upper tail dependence coefficient}. So $X$ and $Y$ are asymptotically independent if and only if their upper tail dependence coefficient is zero.

Sibuya's condition is simple, but the formulation in terms of survival probabilities is inconvenient. The quantiles in~(\ref{q2con}) may be hard to determine since this amounts to computing the inverse of the dfs. There is a simple criterium in terms of sums:
\bpr\label{p2s}
Let $(X_1,X_2)$ have a df with continuous marginals $F_1$ and $F_2$, and let $1-F_i(t_{in})=1/n$. If $s_n:=n\pbb\{X_1+X_2>t_{1n}+t_{2n}\}\to0$ then $X_1$ and $X_2$ are asymptotically independent.
\epr
\bew Sibuya's function $P$ satisfies $nP(1-1/n,1-1/n)=n\pbb\{X_1>t_{1n},X_2>t_{2n}\}\le s_n\to0$. This gives Sibuya's condition~\eqref{q2ai}. \qed

Below we give criteria in terms of continuous curves $\xb(t)=(x_1(t),x_2(t))$, $t\ge0$, for which $F_1(x_1(t))$ and $F_2(x_2(t))$ tend to one as $t\to\nf$. 

\bpr\label{p2sa}
Let $(X_1,X_2)$ have df $F$ with continuous marginals $F_1$ and $F_2$. The components $X_1$ and $X_2$ are asymptotically independent if and only if for any $\e>0$ there exists a continuous curve $\xb(t)$, $t\ge0$, such that $p_i(t)=\pbb\{X>x_i(t)\}$ is positive and vanishes for $t\to\nf$ for $i=1,2$, and such that
\bql{q2xi}
\pbb\{X_1>x_1(t), X_2>x_2(t)\}/p_i(t)<\e \qquad t>t_\e,\ i=1,2.
\eql
\epr
\bew Assume asymptotic independence. One may choose $x_1(t)$ and $x_2(t)$ continuous and increasing such that $p_1(t)=p_2(t)$ for all $t$. Then~\eqref{q2ai} gives~\eqref{q2xi}. Now assume~\eqref{q2xi}. Let $u\in(0,p_0]$ where $p_0=\min\{p_1(t_0),p_2(t_0)\}$. By symmetry we may assume $p_1(t)=u\le p_2(t)$. Then
$$P(1-u,1-u)\le\pbb\{X_1>x_1(t), X_2>x_2(t)\}\le\e p_1(t)=\e u.$$
This holds for all $u\in(0,p_0]$. So Sibuya's condition is satisfied. \qed

Asymptotic independence is preserved under quite severe deformations of the distribution. 
\bcr\label{cai2}
Let $C$ be a convex open cone in $\rbb^d$ and let $f$ and $g$ be probability densities which are positive on $C$, vanish off $C$ and for which the quotients $f/g$ and $g/f$ are bounded on $C$. Let $\ZB=(Z_1,\ldots,Z_d)$ have density $f$ and let $\XB=(X_1,\ldots,X_d)$ have density $g$. If $Z_1$ and $Z_2$ are asymptotically independent, then so are $X_1$ and $X_2$.
\ecr
\bew The inequality $f\le Mg$ implies by integration that the same inequalities hold for the univariate and bivariate marginal densities and tail functions. So $\pbb\{X_1>x_1(t), X_2>x_2(t)\}\le M\pbb\{Z_1>x_1(t), Z_2>x_2(t)\}$. Similarly the inequality $g\le Mf$ gives an inequality with the constant $M$ for the tail probabilities: $\pbb\{Z_i>x_i(t)\}\le M\pbb\{X_i>x_i(t)\}$ for $i=1,2$. Now use~(\ref{q2xi}) for $(Z_1,Z_2)$ to establish the same relation for $(X_1,X_2)$. The extra factor $M^2$ has no effect in the limit. \qed

The concept of asymptotic independence has been refined by looking at the rate at which $P(1-s,1-s)$ vanishes for $s\to0$, or more generally by looking at the behaviour of the survival function $P(1-u,1-v)$ for $u,v\to0$. See~\cite{Resnick2002} or \cite{Ramos2009}. This second order theory is called hidden regular variation. We shall not treat this subject in our paper. Since our interest is in multivariate densities rather than dfs, we give conditions on the density which ensure asymptotic independence of the components. So the assumption of continuous dfs is not restrictive. We shall also assume that the densities are continuous on a convex cone and vanish outside this cone. In our context, this cone will typically be either the whole space $\rbb^d$ or the open positive orthant $(0,\nf)^d$.

\subsection{Asymptotic independence and multivariate extreme value theory}

Asymptotic independence has to do with extremes,   more precisely with bivariate maxima. Knowledge of multivariate extreme value theory  is not indispensable for understanding asymptotic independence, but it will help to better understand asymptotic dependence. 

For multivariate extreme value theory  we have to assume that each of the marginals of the multivariate df $F$ lies in the domain  of attraction of a univariate extreme value limit law (see e.g. Section~0.3 in \cite{Resnick1987} for a definition).  As above let $\ZB_1,\ZB_2,\ldots$ with $\ZB_n=(X_n,Y_n)$ be independent observations from the bivariate df $F$, and write $\ZB^{\lor n}$
for the $n$th coordinatewise partial maximum. For simplicity assume the marginals are equal with the standard Fr\'echet df $F_i(t)=H(t)=e^{-1/t}$ on $(0,\nf)$. The normalized marginal maxima $U_n=X^{\lor n}/n$ and $V_n=Y^{\lor n}/n$ again have the Fr\'echet df $H$ by the scaling property $H^n(nt)=H(t)$. The scaled bivariate maximum $\WB_n=\ZB^{\lor n}/n$ has df $G_n(\wb)=F^n(n\wb)$. Suppose $G_n$ converges weakly to a limit distribution $G$, known as a multivariate extreme value distribution or a \emph{max-stable} distribution. The limit vector $\WB=(U,V)$ lives on $(0,\nf)^2$, and the components $U$ and $V$ have a Fr\'echet law. Asymptotic independence for the df $F$ is equivalent to independence of the components of the limit vector (see e.g. Proposition~5.27 in \cite{Resnick1987} or Theorem~6.2.3 in \cite{Haan2006}).

The convergence $F^n(n\wb)\to G(\wb)$ becomes easier to handle if one takes logarithms. Write 
$G(\wb)=e^{-R(\wb)}$. Since $-\log F(n\wb)$ is asymptotic to $1-F(n\wb)$ for $\wb>\zerob$ we may write the limit relation as
\bql{q2R}
n(1-F(n\wb))\to R(\wb)\qquad \wb\in\XC=[0,\nf)^2\sm\{(0,0)\}.
\eql
The left hand side $R_n(\wb)=n(1-F(n\wb))$ is a df of the measure $\r_n=n\p_n$, where $\p_n$ is the probability distribution of the vector $\ZB/n$. So $\r_n$ is the mean measure of the  sample cloud $N_n=\{\ZB_1/n,\ldots,\ZB_n/n\}$ and $R_n(\wb)=\r_n([\zerob,\wb]^c)$. Pointwise convergence $G_n\to G$ on $[0,\nf)^2$ implies pointwise convergence $R_n\to R$ on $\XC$ and vague convergence $\r_n\to\r$. One can prove weak convergence $\r_n\to\r$ on $\XC\sm[\zerob,\wb]$ for any $\wb\in(0,\nf)^2$. It follows that the scaled sample clouds $N_n$ converge in distribution to a Poisson point process $N$ on $\XC$ weakly on the complement of centered disks.  The Poisson point process $N$ is of interest since it gives an asymptotic description of the large vectors in the sample cloud. Moreover, the coordinatewise maximum $\WB_n$ of the scaled sample cloud $N_n$ converges in distribution to the coordinatewise maximum $\WB$ of $N$ as the sample size goes to infinity.

Sibuya's condition holds precisely if $R(u,v)=1/u+1/v$. (Indeed $n(1-F_i(nt))\to1/t$ for $i=1,2$, and the survival function $\ov F(x,y)=\pbb\{X>x,Y>y\}=P(F_1(x),F_2(y))$ satisfies $n\ov F(n\wb)\to0$ by Sibuya's condition and Proposition~\ref{p2ab}.) In this case the measure $\r$ has density zero on the open quadrant $(0,\nf)^2$  by differentiation (see also Proposition~5.24 in \cite{Resnick1987}). So $\r$ is the sum of two  measures, one on the positive horizontal and one the positive vertical axis, both with density $1/t^2$. The point process $N$ thus is the superposition of two Poisson point processes on these halfaxes, and the two point processes are independent since the halfaxes are disjoint. This yields a simple description of the behaviour of large sample clouds: Asymptotically there is no relation between very large observations in the horizontal direction, and very large observations in the vertical direction. 

So far we have looked at asymptotic independence for bivariate distributions. Unlike independence, for a multivariate df asymptotic independence holds if it holds for the bivariate marginals; see e.g. Remark~6.2.5 in \cite{Haan2006}. Here is the argument. For simplicity assume $F$ has standard Fr\'echet marginals, and the multivariate maxima converge in distribution. The scaled sample clouds $N_n=\{\ZB_1/n,\ldots,\ZB_n/n\}$ then converge in distribution to a Poisson point process $N$ on $\XC=[0,\nf)^d\sm\{\zerob\}$ with mean measure $\r$ weakly on the complement of any centered ball. If Sibuya's condition holds for the bivariate marginals, then the projection $\r_{ij}$ of $\r$ on the $x_i,x_j$-plane will live on the two positive halfaxes in this plane. It follows that $\r$ lives on the set of points in $[0,\nf)^d\sm\{\zerob\}$ which have at most one positive coordinate, the union of the $d$ positive halfaxes. The restrictions of $N$ to these halfaxes are independent (since the $d$ positive halfaxes are disjoint). Hence the multivariate extreme value limit vector has independent components. 

\subsection{Asymptotic dependence}

Asymptotic dependence is an ambiguous term. Logically it means the absence of asymptotic independence. We shall usually interpret it in a more constructive manner to mean the existence of a max-stable limit law $G=e^{-R}$. This in turn implies convergence in distribution of the normalized sample clouds to a Poisson point process $N$, whose mean measure $\r$ has df $R=-\log G$, as sketched above. The tail dependence coefficient $\l_U$ in~\eqref{q2con} gives only a very restricted view of extremal dependence. As pointed out in Section~8.2 in \cite{Mikosch2006}, a positive coefficient does not imply that the underlying distribution belongs to the maximum domain of attraction of some extreme value limit law. The exponent measure $\r$ provides a much more informative description of asymptotic dependence. Points of $N$ which do not lie on one of the axes denote a simultaneous occurrence of very large values in two or more coordinates in the corresponding point of the normalized sample. Thus the point process $N$ (or equivalently $\r$) gives a complete description of how the extreme upper order statistics in the different coordinates are linked.

\section{Densities and level sets}

The aim of this section is to give conditions in terms of the density which will guarantee asymptotic independence. Two aspects of densities play an important role in our analysis: the shape of the level sets and the tail behaviour. We consider densities  which are completely specified by just these two quantities - a shape for the level sets, and a decreasing function governing the rate of decay of the tails. In many cases of practical interest the shape is a bounded open convex set, containing the origin. The density then is continuous if and only if the decreasing function is. We shall also consider such unimodal densities which vanish outside the positive orthant.

\subsection{Homothetic densities}

Homothetic densities are densities whose level sets are scaled copies of a given open set $D$: 
\bql{erc}
\{f>c\}:=\{\xb\in\rbb^d\mid f(\xb)>c\}=r_cD\qquad 0<c<c_0:=\sup f.
\eql
We assume that $D\ss\rbb^d$ is a bounded open star-shaped set. A \emph{star-shaped set} has the property that with any point $\xb$ it contains all points $r\xb$ for $0<r<1$. Assume the set $D$ contains the origin. If each ray intersects the boundary $\prl D$ in one point, then the set $D$ can be represented using a \emph{gauge function} $n_D:\rbb^d\to[0,\nf)$ that satisfies (\rm{i}) $n_D(t\xb)=tn_D(\xb)$ for $t>0$, $\xb\in\rbb^d$ (homogeneity property), and (\rm{ii}) $D=\{\xb\in\rbb^d\mid n_D(\xb)<1\}$. The conditions on the set $D$ and the continuity of $f$ ensure that the gauge function is continuous. If the set $D$ is convex, then so is the gauge function. If in addition the set is symmetric, $-D=D$, then the gauge function is a norm on $\rbb^d$, and the set $D$ is the open unit ball in this norm. For any bounded open star-shaped set $D$ the sets $nD$, $n>0$, form an increasing family. Their union $D_\nf$ is an open cone. It is a proper cone if the origin is a boundary point of $D$. One may then still define the gauge function $n_D$. This now is a function on the open cone $D_\nf$. It is continuous on the cone if each ray in the cone intersects the boundary of $D$ in a single point, see Proposition~\ref{p3c} below.

Gauge functions allow us to give analytic expressions for homothetic densities with the same ease with which one handles spherically symmetric densities. Continuous homothetic densities have the form $f(\xb)=f_0(n_D(\xb))$ for a decreasing continuous function $f_0$ on $(0,\nf)$. This formula holds for all $\xb$ for which the ray through $\xb$ intersects the set $D$. 

\bdf Let $\DC_d$ denote the class of all bounded open star-shaped sets $D\ss\rbb^d$, for which the cone $D_\nf=\bigcup_nnD$ is convex, and for which the gauge function $n_D$ is continuous on this cone. A density $f$ on $\rbb^d$ belongs to the class $\HC(D)$ if the \emph{shape} $D$ belongs to the class $\DC_d$, and if $f$ is of the form $f(\xb)=f_0(n_D(\xb))$,
where the \emph{density generator} $f_0:[0,\nf)\to[0,\nf)$ is decreasing, positive and continuous. We set $f\equiv0$ outside the cone $D_\nf$ on which the gauge function is defined.
\edf

It is apparent from the above definition that densities in $\HC$ are \emph{(star) unimodal} since all the level sets $\{f>c\}$, $c\in(0,c_0)$ are star-shaped; cf. Section~2.2 in \cite{Dharmadhikari1988}. Typical examples of densities in $\HC$ to keep in mind are the multivariate centered normal densities and, more generally, elliptically symmetric densities, discussed for example in \cite{Fang1990}. See also Example~\ref{e60}. In certain applications, elliptical symmetry may be too restrictive. Densities in $\HC$ give the flexibility to model directional irregularities present in the data clouds, and to handle distributions on the positive orthant. The regularity conditions exclude pathological sets; see Example~\ref{e30}.

Before we proceed by looking at properties of these densities, let us review some related classes of models proposed in the literature. The $\lp$-spherical densities (\cite{Osiewalski1993}) extend the class of spherical densities by allowing level sets to be balls in $\lp$-norm for any $p\ge1$. A further generalization is given by so-called $v$-spherical densities (\cite{Fernandez1995}), where the scale function $v$ plays the same role as the gauge function $n_D$ defined above. In fact, our class $\HC(D)$ is a subclass of the $v$-spherical densities in that we restrict level sets to be bounded and star-shaped, and $f$ to be continuous on $D_\nf$. In a recent paper \cite{Arnold2008}, the authors advocate to study densities in terms of their contours.

For a given shape $D$, what conditions does the density generator $f_0$ have to satisfy in order that the function $f(\cdot)=f_0(n_D(\cdot))$ is a probability density on $\rbb^d$? By regarding the set below the graph of $f$ as a pile of thin $D$-shaped slices we obtain the following partial integration result:
\bql{qpd}\pbb\{\XB\in tD\}=\int_{tD}f_0(n_D(\xb))d\xb=f_0(t)|tD|+|D|\int_0^ts^d|df_0(s)|=\int_0^tf_0(s)d|sD|,\quad t>0.\eql
Observing that $|sD|=s^d|D|$ and letting $t$ tend to infinity, we obtain the condition (cf. Equation~(5) in \cite{Fernandez1995}): $1=d|D|\int_0^\nf s^{d-1}f_0(s)ds$.

The class of unimodal densities introduced above is invariant under linear transformations. If the vector $\XB$ has density $f\in\HC(D)$ then the vector $\YB=A\XB$ has density $g\in\HC(E)$, where $E$ is the image of the star-shaped (convex) set $D$ under the linear transformation $A$, and hence is also star-shaped (convex). A nice illustration of this invariance is the extension of spherical distributions to elliptical ones.

For densities in $\HC$ the distinction between light and heavy  tails is crucial for asymptotic independence. It is determined by the behaviour of the generator $f_0$ at infinity. Let us recall the definitions of regular and rapid variation. 
\bdf
A measurable function $h:(0,\nf)\to(0,\nf)$ is \emph{regularly varying} at $\nf$ with exponent $\q$, if for $x>0$
\bql{erev}\lim_{t\to\nf}h(tx)/h(t)=x^\q\qquad\q\in\rbb;\eql
if $\q=0$, then $h$ is called \emph{slowly varying};
$h$ is \emph{rapidly varying} at $\nf$ if 
\bql{erav}\lim_{t\to\nf}h(tx)/h(t)=\bcs\nf, & 0<x<1\\ 0,& x>1.\ecs\eql  
\edf

If the generator of a density $f=f_0(n_D)$ varies rapidly, the density has light tails. In our terminology $f$ has \emph{light tails} precisely if $f_0$ is continuous, positive, and strictly decreasing on $[0,\nf)$, and if there is a strictly increasing sequence $r_n>0$ such that $r_{n+1}\sim r_n$ and $f_0(r_n)/f_0(r_{n+1})\to\nf$.

Rapid variation of $f_0$ allows us to give strong inequalities for the measure $\m$ with density $f=f_0(n_D)$. 
\bpr\label{pAC1}
Let $\m$ have density $f\in\HC(D)$ with a rapidly varying density generator $f_0$.
\bit
\item[{\rm(i)}] For any $\e>0$
\bql{qACm}\m(rD^c)<<\m(rD\sm e^{-\e}rD)\qquad r\to\nf.\eql
\item[{\rm(ii)}] For any non-empty open set $U\ss D$
\bql{qACinU} \m(rD^c)<<\m(rU)\qquad r\to\nf.\eql
\eit\epr
\bew Rapid variation implies that $f_0(e^\d t)<f_0(t)/M$ eventually for $t\to\nf$, and hence the rings $R_n=e^{n\d+\d}tD\sm e^{n\d}tD$, $n\ge0$, have measure $\m(R_{n+2})\le e^{2d\d}\m(R_n)/M$. The rings $R_n$ are disjoint, and their union is the complement of $tD$. On summing the odd and the even terms we find with $\h=e^{2d\d}/(M-e^{2d\d})$:
$$\m(e^{2\d}tD^c)=\sum_{n=1}^\nf(\m(R_{2n})+\m(R_{2n+1}))\le\h(\m(R_0)+\m(R_1))=\h\m(e^{2\d}tD\sm tD).$$
This gives~\eqref{qACm}. The integral over a thin ring is much larger than the integral over the set outside the ring. Formally, take $\d=\e/2$ and $r=t-\e$ to obtain~\eqref{qACm}. To prove~\eqref{qACinU} take a non-empty open subset $U_0\ss U$ whose closure lies in $D$. Then $e^{2\e}U_0\ss D$ if $\e$ is small, and hence by rapid variation of $f_0$ as above the infimum of $f$ over $rU_0$ is much larger than the supremum of $f$ over the ring $R=rD\sm e^{-\e}rD$. 
Since $|R|/|rU_0|$ is a constant we conclude that $\m(R)<<\m(rU_0)$, and~(\ref{qACinU}) follows from~(\ref{qACm}).\qed


We now give some extra details on star-shaped sets. 
\bpr\label{p3c}
Let $D$ be a bounded open star-shaped set. Suppose for each non-zero vector $\xb\in D$ there is one positive real $r_0$ such that $r_0\xb$ lies on the boundary of $D$. Then the gauge function $n_D$ is continuous on the cone $D_\nf=\bigcup_n nD$. 
\epr
\bew Let $B$ denote an open unit ball. If $D$ contains the origin, it contains a ball $\e B$, and by homogeneity $n_D<\d$ on $\d\e B$. So the gauge function then is continuous in the origin. Now suppose it is discontinuous in a point $\pb$ outside the origin. We may assume that $n_D(\pb)=1$ and that there is a sequence $\pb_n\in D_\nf$ such that $\pb_n\to\pb$ and $n_D(\pb_n)\to c\ne1$. Let $r\in(0,1)$ be close to $1$. Then $n_D(r\pb_n)\to rc$, and $r\pb\in D$ implies $r\pb_n\in D$ eventually, hence $rc<1$ and so $c<1$ since $r$ is arbitrary. This implies $\pb_n\in D$ eventually, and this also holds for $s\pb_n$ for $1<s<1/c$. Hence $s\pb$ is a boundary point. But so is $\pb$. This contradicts our assumption on the boundary of $D$. \qed

In general densities in $\HC(D)$ are not closed under projection, even if $D$ is convex. If a random vector $(X_1,\ldots,X_d)$ has density $f\in\HC(D)$, then the density of $(X_1,\ldots,X_{d-1})$ need not be homothetic, and the univariate marginals need not even be unimodal; see Example~\ref{e40}. There are some exceptions. Projections of spherical densities are spherical, and if the level sets are balls in $\lp$ for some $p\in[1,\nf]$ then this also holds for projections along the coordinate axes. The class $\DC_d$ itself is closed under projection. If $D$ is a bounded open convex set containing the origin, then so is the vertical projection $E$ of $D$ onto the horizontal hyperplane; if the origin is a boundary point of $D$ it may be an interior point of $E$.

\bpr\label{p3p}
Suppose $D\in\DC_d$. Let $E$ be the vertical projection of $D$ onto the horizontal hyperplane. Then $E\in\DC_{d-1}$.
\epr
\bew Write $\zb=(\xb,y)$ to distinguish the horizontal and vertical part of the vector $\zb$. It is clear that $E$ is a bounded open star-shaped set. Moreover $E_\nf$ is the projection of the cone $D_\nf$, and hence an open convex cone. We have to prove continuity of the gauge function $n_E$ on $E_\nf$. If $n_E$ is not continuous, there exists a vector $\xb$ and a sequence $\xb_n\to\xb$ such that $n_E(\xb)=1$ and $n_E(\xb_n)<c_0<1$ because $E$ is open. Hence $n_D(\xb,y)\ge1$ for all $y$ and there exist $y_n$ such that $n_D(\xb_n,y_n)<c_0$. Since $D$ is bounded, the sequence $(y_n)$ is bounded, and we may assume that it converges to some element $y_0$. Then $(\xb_n,y_n)\to(\xb,y_0)$, and continuity of $n_D$ implies $n_D(\xb,y_0)\le c_0<1$. Contradiction. \qed

\subsection{Densities whose level sets are asymptotically star-shaped}

In this section we relax the condition that all level sets have the same shape $D\in\DC_d$ to the condition that the level sets, properly scaled, converge to a set $D\in\DC_d$. We restrict attention to the light-tailed setting. 

\bdf\label{dsae}
Let $D\in\DC_d$. A positive probability density $f$ on $\rbb^d$ belongs to the set $\AC(D)$ if there exist sequences $c_n>0$ and $r_n\to\nf$ with $c_{n+1}/c_n\to0$ and $r_{n+1}\sim r_n$ such that for any $\e>0$ eventually
\bql{qag1} e^{-\e}r_nD\ss\{f>c_n\}\ss e^\e r_nD,\qquad n\ge n_0.\eql
We write $\{f>c_n\}/r_n\to D$. A continuous positive function $\tilde f$ is \emph{shape equivalent} to $f$ if its level sets satisfy \eqref{qag1}. \qend
\edf

The sequences $r_n$ and $c_n$ determine a set of continuous decreasing functions $\h$ which satisfy $\h(c_n)\sim r_n$. All these functions $\h(c)$ vary slowly for $c\to0^+$ by the assumption that $c_{n+1}/c_n\to0$ and $r_{n+1}/r_n\to1$. It is this set of slowly varying functions rather than the particular sequences $c_n$ and $r_n$ which are of interest.

\bpr\label{pach} 
If the slowly varying function $\h$ above is strictly decreasing, defined on $(0,c_\h]$ for some $c_\h>0$, and vanishes in $c_\h$, then the inverse function $f_0=\h\ginv$ is a continuous positive strictly decreasing function on $[0,\nf)$ which varies rapidly in $\nf$, and $f=f_0(n_D)$ is shape equivalent to $g$, and for all $\e>0$ 
\bql{qag2}e^{-\e}\{f>c\}\ss\{g>c\}\ss e^\e\{f>c\}\qquad 0<c<c_\e.\eql
\epr
\bew Rapid variation of the inverse function $f_0$ follows from Theorem~2.4.7(i) in \cite{Bingham1987}. Shape equivalence holds since $\{f>c\}=\h(c)D$ and $\h(c_n)\sim r_n$ by assumption.
\qed

Since $r_n\to\nf$ one may take a strictly increasing subsequence such that the asymptotic equality $r_{n+1}\sim r_n$ remains valid. Take $c_0>c_1$. Any continuous strictly decreasing function $\h$ on $(0,c_0]$ which vanishes in $c_0$ and has the value $r_n$ in $c_n$ satisfies the conditions of the proposition above. So there are many continuous strictly decreasing functions $f_0$ on $[0,\nf)$ which vary rapidly in $\nf$ such that $f_0(n_D)$ is shape equivalent to $g$. 

For functions in $\AC(D)$ the inequalities in Proposition~\ref{pAC1} also hold. Moreover, they have the nice property that the $d$ marginals $g_i$ of $g\in\AC(D)$ will lie in $\AC(D_i)$, where $D_i$ is the projection of $D$ on the $i$th coordinate.

\blm\label{pAC2} Let $g\in\AC(D)$. Let $\r_r$ for $r\ge1$ denote the measure with density $\ub\mapsto g(r\ub)/g(r\qb)$ for a fixed non-zero vector $\qb$. Then for any open set $U$ which intersects $D$
\bql{qACinB} \r_r(D^c)<<\r_r(U).\eql
\elm
\bew Let $f=f_0(n_D)$ be shape equivalent to $g$ and define $\m_r$ to have density $f_r(\ub)=f(r\ub)/f_0(r)$.  It suffices to prove~\eqref{qACinB} for the measures $\r_r$ with density $g_r(\ub)=g(r\ub)/f_0(r)$. Let $V$ be a non-empty open subset of $U$ which lies in the complement of $\d_0D$ for some $\d_0>0$, and whose closure lies in $D$. For $r>r_0/\d$ we may apply the pointwise inequality~\eqref{qACfg} below and conclude that $\r_r(U)\ge e^{-\e d}\m_r(e^\e V)$, where we choose $\e>0$ so small that $e^{3\e}V\ss D$. Then, using \eqref{qACinU}, $\m_r(e^{-\e}D^c)<<\m_r(e^\e V)$ and another application of the pointwise inequality gives~\eqref{qACinB}. \qed

\bpr\label{p3Y}
Suppose $\XB$ has density $f\in\AC(D)$. Let $Y=\x(\XB)$ be a non-zero linear combination of the components of $\XB$. Then $Y$ has a density $g\in\AC(J)$ where $J$ is the open interval $\x(D)$.
\epr
\bew We may assume that $\x$ is the vertical coordinate. The condition $f(e^\e\xb)<f(\xb)/M$ for $\|\xb\|\ge r_0$ implies by integration over horizontal hyperplanes $g(e^\e y)<e^{(d-1)\e}g(y)/M$ for $|y|\ge r_0$. This gives rapid variation. The thin tails of $f$ ensure that $g$ is continuous. Let $\x(D)=(a,b)$. We claim that the average of $g$ over $r(a,a+\e a)$ and over $r(b-\e b,b)$, say $a_1$ and $b_1$, is much larger than over both $r(a-\e a,a)$ and $r(b,b+\e b)$, say $a_2$ and $b_2$, as $r\to\nf$. It suffices to show that $\pbb\{ra<Y<ra+\e ra\}$ and $\pbb\{rb-r\e b<Y<rb\}$ are much larger than  $\pbb\{\XB\in rD^c\}$. This follows from~\eqref{qACinB} since the horizontal slice $\{ra<x_d<ra+\e ra\}$ contains the open set $U=rD\cap\{x_d<ra+\e ra\}$, and $\pbb\{\XB\in U\}$ is much larger than $\pbb\{\XB\in rD^c\}$. A similar argument holds for the strip $\{rb-\e rb<x_d<rb\}$. Let $c$ lie between $a_1\land a_2$ and $b_1\lor b_2$. Then $g>c$ on $[ra-2\e ra,rb+2\e rb]$ and $g<c$ holds off $[ra-2\e ra,rb+2\e rb]$. Hence $e^{-3\e}rJ\ss\{g>c\}\ss e^{3\e}rJ$ holds for small $\e$ for sufficiently large $r$. \qed

It is not known whether a similar result holds for the bivariate marginals of $f$.


One can also define shape equivalence for functions with heavy tails. In that case shape equivalence is the same as asymptotic equality in infinity, see Proposition~\ref{pau1}. For light tails the behaviour of the quotient of two shape equivalent functions may be very erratic. Example~\ref{eAAC} in the Appendix exhibits some functions which are shape equivalent to the bivariate Gaussian density.

The level set inequalities in \eqref{qag2} imply the pointwise inequalities:
\bql{qACfg} f(e^\e\xb)\le g(\xb)\le f(e^{-\e}\xb)\qquad\xb\in r_\e D^c.\eql 
This shows that a function $g=qf$ is shape equivalent to $f=f_0(n_D)$ for $q=e^\j$ and $f_0=e^{-\f}$ if for every $\e>0$
\bql{qACq} \f(r)-\f(e^\e r)\le\j(\xb)\le\f(r)-\f(e^{-\e}r)\qquad\xb\in\prl rD,\quad r\ge r_\e.\eql
Rapid variation implies that the left side goes to $-\nf$ and the right to $\nf$. Hence we find:

\bpr $g\in\AC(D)$ is shape equivalent to $\tilde g=qg$ if $q$ is positive and $\log q$ is bounded. 
\epr

For Weibull-like functions of the form $f=f_0(n_D)$ where $f_0=e^{-\f}$ for a continuous strictly increasing function $\f$ which varies regularly with positive exponent there is a simple alternative description of shape equivalence. 

\blm\label{prve}
Suppose $f_0=e^{-\f}$ where $\f$ is continuous, strictly increasing and varies regularly in $\nf$ with exponent $\q>0$. Then $f_0$ varies rapidly in $\nf$ and $g=e^{-\g}\in\AC(D)$ is shape equivalent to $f_0(n_D)$ if and only if $\g(\xb_n)\sim \f(n_D(\xb_n))$ holds for all sequences $\|\xb_n\|\to\nf$. 
\elm
\bew Regular variation of $\f$ gives $\f(e^\e r)-\f(r)\sim (e^{\q\e}-1) \f(r)$ and $\f(r)-\f(e^{-\e}r)\sim(1-e^{-\q\e})\f(r)$ for $r\to\nf$. The claim then follows from \eqref{qACq} with $\g=\f-\j$.
\qed

\bpr\label{prd1}
Suppose $g=e^{-\g}$ is a continuous positive density on $\rbb^d$. Let there exist a function $\n$ on $\rbb^d$ which is positive outside a bounded set and a non-zero vector $\qb$ such that 
\bql{emrv1}\g(t_n\ub_n)/\g(t_n\qb)\to\n(\ub)\qquad t_n\to\nf,\ \ub_n\to\ub,\ \ub\in\rbb^d.\eql
Then there exists a set $D\in\DC_d$ containing the origin, a positive constant $\q$ such that $\n=n_D^\q$, and a continuous strictly increasing function $\f$ on $[0,\nf)$ which varies regularly with exponent $\q$ in $\nf$ such that $\g$ is asymptotic to $\f(n_D)$ in $\nf$, and $g$ is shape equivalent to $f_0(n_D)$ where $f_0=e^{-\f}$. 
\epr
\bew Set $a(t)=\g(t\qb)$. Then $a$ varies regularly in $\nf$ with exponent $\q\ge0$. For any unit vector $\w$ the function $t\mapsto \g(t\w)$ varies regularly in $\nf$ since $\g(t\w_n)/a(t)$ has a positive limit $\n(\w)$ for $\w_n\to\w$ and $t\to\nf$. Uniform convergence on compact sets implies that $\n$ is continuous on the unit sphere, and $a(st)/a(t)\to s^\q$ implies $\g(st\w)/a(t)\to\n(\w)s^\q$. This proves that $\n(r\w)=r^\q\n(\w)$, and $\g(r\w)\sim a(r)\n(\w)\asymp a(r)$ for $\|r\w\|\to\nf$ since~\eqref{emrv1} implies uniform convergence on compact sets, and $\n$ is continuous and hence bounded on compact sets. Hence $a(r)\to\nf$ and~\eqref{emrv1} with $\ub_n=\zerob$ gives $\n(\zerob)=0$ and continuity of $\n$ implies $\q>0$. Set $D=\{\n<1\}$. Then $\n=n_D^\q$. Let $\pb\in\prl D$. Then $t\mapsto \g(t\pb)$ varies regularly with exponent $\q$, and we may choose $\f$ positive, continuous and strictly increasing and asymptotic to this function for $t\to\nf$. Then $\g(r\w)\sim a(r)\n(\w)$ gives $\g(\wb)\sim \f(n_D(\wb))$ for $\|\wb\|\to\nf$. Lemma~\ref{prve} then shows that $g$ is shape equivalent to $e^{-\f(n_D)}$ and since $e^{-\f}$ varies rapidly if $\f$ varies regularly with exponent $\q>0$ it follows that $g\in\AC(D)$. \qed

\brm Condition~\eqref{emrv1} is related to multivariate regular variation; see e.g. Section~5.4.2 in \cite{Resnick1987}.
\erm

\subsection{Criteria for asymptotic independence}

A random vector $\XB$ with a spherically symmetric density $f(\xb)=f_0(\|\xb\|_2)$ will have asymptotically independent components if the generator $f_0$ varies rapidly, and it will have asymptotically dependent components if $f_0$ varies regularly (see e.g. Theorem~4.3 in \cite{Hult2002}, or Proposition~3.2 in \cite{Hashorva2005a}). These results remain valid for $f\in\AC(B)$ where $B$ denotes the open unit ball. One may replace $B$ by certain bounded open star-shaped sets $D$, as will be shown in Theorem~\ref{thmai1} and \ref{thmad1} below. 

With any bounded open set $D$ one may associate the open intervals $D_i=(a_i,b_i)$, $i=1,\ldots,d$, obtained by projecting $D$ onto the $i$th coordinate. Then $(\ab,\bb)$ is the smallest open box containing $D$, and $\bb=(b_1,\ldots,b_d)=\sup D$ is the coordinatewise supremum of all points in $D$ (and $\ab=(a_1,\ldots,a_d)=\inf D$). 

Asymptotic independence depends on the bivariate marginal distributions. Hence we introduce the projections of $D$ on the $(x_i,x_j)$ coordinates. For $1\le i<j\le d$ we denote by $D_{ij}\ss\rbb^2$ the projection of $D$ onto the two-dimensional space spanned by the unit base vectors $\eb_i$ and $\eb_j$. The sets $D_{ij}$ lie in $\DC_2$ by Proposition~\ref{p3p}, and $D_{ij}$ fits exactly into the rectangle $(a_i,b_i)\times(a_j,b_j)$. The cone generated by $D_{ij}$ is the projection of the cone $D_\nf$ generated by $D$. 

\bdf\label{dblunt} The set $D\in\DC_2$ is \emph{blunt} if the point $(b_1,b_2)=\sup D$ does not lie in the closure of $D$. \edf 
A bounded open convex set $D$ in $\rbb^d$ is \emph{smooth} in the boundary point $\pb$ if there is a unique hyperplane which contains $\pb$ but which does not intersect $D$, the tangent plane to $D$ at $\pb$. For a planar set $D$ this means that $\pb$ is not a vertex. If the convex hull of $D$ is smooth in all points then all bivariate projections $D_{ij}$ are blunt. 

We can now state our main results for asymptotic independence in terms of densities.

\bth\label{thmai3}
If $\XB$ has a light-tailed homothetic density $f\in\HC(D)$, and $D$ is convex with a smooth boundary, then for any two distinct unit vector $\ab$ and $\bb$ the random variables $\ab^T\XB$ and $\bb^T\XB$ are asymptotically independent. The result remains valid if the density of $\XB$ is in $\AC(D)$.
\eth
\bew First assume $\ab$ and $\bb$ are linearly independent. Introduce new coordinates such that $\ab$ and $\bb$ become the first two base vectors $\eb_1$ and $\eb_2$. It suffices to check that the assumption holds for vertical tangent planes, hyperplanes which project onto a line in the two-dimensional $x_1,x_2$-plane. The characterization of $D$ is geometrical and so it is preserved under linear transformations. The projection $D_{12}$ in the new coordinates is also convex and smooth. Hence $D_{12}$ is blunt, and we may apply Theorem~\ref{thmai1} below. If $\bb=-\ab$ the bivariate distribution lies on the counterdiagonal,  $y=-x$, and asymptotic independence is trivial by applying Proposition~\ref{p2sa} with $\xb(t)=(t,t)$ for $t\ge0$. \qed

\bth\label{thmai1}
Suppose $\XB$ has density $g\in\AC(D)$. If the bivariate projection $D_{12}$ is blunt then $X_1$ and $X_2$ are asymptotically independent. 
\eth
\bew There is a simple analytic argument. Let $\sup D_{12}=(b_1,b_2)$. The sum $Y=X_1+X_2$ has density $g\in\AC(J)$ by Proposition~\ref{p3Y} where $J$ has upper endpoint $b=\sup\{x+y\mid (x,y)\in D_{12}\}$ and $b<b_1+b_2$ by bluntness. Now use Sibuya's condition for sums in Proposition~\ref{p2s}. We give a more probabilistic proof in Section~\ref{s:sc}.\qed

\subsection{Criteria for asymptotic dependence}

We now give the counterpart to Theorem~\ref{thmai1} for heavy-tailed densities. We are interested in the case where the partial maxima go to infinity in all coordinates. If the shape $D$ lies in the negative orthant then the coordinatewise maxima converge to $\zerob$, if it lies in a negative coordinate halfspace $\{x_i<0\}$ all partial maxima will lie in this halfspace too. We exclude these cases in the theorem below.

\bth\label{thmad1}
Suppose $\XB$ has density $f\sim f_0(n_D)$ with $D\in\DC_d$ and $f_0$ a continuous strictly decreasing positive function on $[0,\nf)$ which varies regularly in $\nf$ with exponent $-(\l+d)$ for some $\l>0$. 
Assume for each coordinate $i\in\{1,\ldots,d\}$ the set $D$ contains a point whose $i$th component is positive. The components of $\XB$ are asymptotically dependent unless $D$ is contained in the set $S$ of points with at most one positive coordinate, a union of $d+1$ orthants. The partial maxima $\XB^{\lor n}$ may be scaled to converge in law to a vector $\WB$ whose components have df $\pbb\{W_i\le t\}=e^{-(a_i/t)^\l}$ for  positive constants $a_1,\ldots,a_d$. The exponent measure $\r^+$ of $\WB$ is the image under the map 
$\xb\mapsto\xb^+=(x_1\lor0,\ldots, x_d\lor0)$ of the excess measure $\r$ with intensity $c/n_D^{\l+d}$. One may take $c=1$ by a suitable choice of the scaling constants for the maxima.
\eth
\bew Choose a point $\wb_0$ in the cone $D_\nf$ on the boundary of $D$, and for $r\ge1$ set
$h_r(\wb)=f(r\wb)/f(r\wb_0)=f_0(rn_D(\wb))/f_0(r)$. By regular variation for $\wb\in D_\nf$, $\wb\ne\zerob$:
$$h_{r_n}(\wb_n)\to h(\wb)=n_D(\wb)^{-(\l+d)}\qquad\wb_n\to\wb,\;r_n\to\nf.$$
Convergence $h_r\to h$ holds uniformly on the intersection of the cone $D_\nf$ with any ring $rB\sm\e B$, where $B$ denotes the open unit ball. (The function $1/n_D$ is bounded on such sets.) By Potter's theorem (Theorem~1.5.6 in \cite{Bingham1987}) for any $\e>0$ there exists $r_\e$ such that
\bql{q3rvb}
f_0(rs)/f_0(r)\le2s^\e/s^{\l+d}\qquad r\ge r_\e,\; s\ge1.
\eql
This yields an integrable majorant for the convergence $h_r\to h$ on $D_\nf\sm B$. Lebesgue's dominated convergence theorem implies that $h_r\to h$ in $\LB^1$ on $D_\nf\sm B$, and because of the uniform convergence above also on $D_\nf\sm\e B$ for any $\e>0$. 

Let $\r(r)$ be the finite measure with density $h_r$, and choose $r_n$ so that $\r(r_n)$ has mass $n$. Then $\r(r_n)$ is the mean measure of the scaled sample cloud $N_n=\{\XB_1/r_n,\ldots,\XB_n/r_n\}$, and $\r(r_n)\to\r$ weakly on $D_\nf\sm\e B$ implies $N_n\imp N$ weakly on $D_\nf\sm\e B$ where $N$ is the Poisson point process on $D_\nf$ with intensity $h$. This tells us that the maxima converge. The measure $\r$ is an excess measure on $\rbb^d\sm\{\zerob\}$:
$$\r(rA)=\r(A)/r^\l\qquad r>0,\ A{\rm\ a\ Borel\ set\ in\ }\rbb^d\sm\{\zerob\}.$$
For halfspaces $A=\{x_i\ge1\}$ or $A=\{x_i\le-1\}$ this relation also holds and implies that the marginals $\r_i$ of $\r$ satisfy the same relation, and hence there exist non-negative constants $c_i^\pm$ such that
\bql{q3ci}\r_i[r,\nf)=c_i^+/r^\l\qquad\r_i(-\nf,-r]=c_i^-/r^\l\qquad r>0.\eql
For the coordinatewise maxima of heavy-tailed distributions it is convenient to work on the non-negative orthant, and replace the vector $\XB$ by $\XB^+$, where we use the continuous map $(x_1,\ldots,x_d)=\xb\mapsto\xb^+=(x_1\lor0,\ldots,x_d\lor0)$. We shall write $N_n^+$ and $N^+$ for the images of $N_n$ and $N$ under this map, and denote the mean measures by $\r^+_n$ and $\r^+$. By assumption each component has a positive probability of being positive. Hence $\WB_n:=\max N^+_n=\max N_n\imp\max N=\max N^+=:\WB$. Write $R(\wb)=\r((-\nfb,\wb]^c)=\r^+([\zerob,\wb]^c)$, $\zerob\ne\wb\ge\zerob$. Then the limit distribution $H$ of the coordinatewise maxima is
$$\pbb\{\WB\le\wb\}=\pbb\{N([\zerob,\wb]^c)=0\}=e^{-\r^+([\zerob,\wb]^c)}=e^{-R(\wb)}.$$
Using~\eqref{q3ci} the same argument gives the marginals $H_i(t)=e^{-\r_i[t,\nf)}=e^{-c_i^+/t^\l}$ with $a_i^\l=c_i^+$. Since $\xb^+$ lies on a coordinate axis precisely if $\xb\in S$, the exponent measure $\r^+$ lives on the positive halfaxes if and only if $D\ss S$.\qed


\section{Sample clouds}\label{s:sc}

In this section we look at the asymptotic behaviour of clouds of independent observations from a given light-tailed distribution as the number of data points in the sample approaches infinity. In particular we are interested in the limiting shape of these sample clouds under suitable scaling. Remark that sample clouds can be viewed as finite point processes with a fixed number of points. The motivation for looking at sample clouds is threefold. First of all, there is a relation between the asymptotic shape of the level sets of the underlying light-tailed density and the limit set onto which corresponding scaled sample clouds converge; see Proposition~\ref{pssc} below. Secondly, the point process approach will yield an intuitive proof of our main results. Finally, for sample clouds projection on the horizontal hyperplane is simple: just delete the last coordinate for each sample point.

As before we consider a sequence of i.i.d. random vectors $\{\XB_1,\XB_2,\ldots\}$ from a continuous distribution on $\rbb^d$ with density $f$. Let $N_n=\{\XB_1/s_n,\ldots,\XB_n/s_n\}$ denote a scaled $n$-point random sample (or \emph{sample cloud}) with the scaling constant $s_n>0$, $s_n\to\nf$ for $n\to\nf$. Alternatively, for any Borel set $A\ss\rbb^d$
$$N_n(A)=\sum_{i=1}^n\oneb_A(\XB_i/s_n).$$ The mean measure of $N_n$ is given by $\r_n=n\p_n$, where $\p_n$ is the distribution of the scaled vector $\XB_1/s_n$. The intensity of the $n$-point point process $N_n$ is $h_n(\ub)=n s_n^df(s_n\ub)$ for $\ub\in\rbb^d$.

Under a suitable choice of $s_n$, the scaled observations $\XB_i/s_n$ from a density $f\in\AC(D)$ with high probability will fill out the closure of the shape set $D$ in the sense of the following definition. 

\bdf\label{dconv2} Let $E$ be a compact set in $\rbb^d$ and $\m_n$ finite measures. We say that the measures $\m_n$ \emph{converge onto} $E$ if $\m_n(\pb+\e B)\to\nf$ for any $\e$-ball centered in a point $\pb\in E$, and if $\m_n(U^c)\to0$ for all open sets $U$ containing $E$. The finite point processes $N_n$ \emph{converge onto} $E$ if $\pbb\{N_n(U^c)>0\}\to0$ for open sets $U$ containing $E$, and if
$\pbb\{N_n(\pb+\e B)>m\}\to1$, $m\ge1$, $\e>0$, $\pb\in E$.
\edf We call the set $E$ in the definition above the \emph{limit set}. In fact, the limit set, if it exists, is always star-shaped (see Proposition~4.1 in \cite{Kinoshita1991}).
The following simple criterion is useful for checking convergence (in probability) of scaled sample clouds (see \cite{Balkema2009} for a proof).
\bpr\label{pconv2}
If $N_n$ is an $n$-point sample cloud from a probability distribution $\p_n$ on $\rbb^d$, then $N_n$ converges onto $E$ if the mean measures $\m_n=n\p_n$ converge onto $E$.
\epr

The next theorem gives a sufficient condition for asymptotic independence of a distribution in terms of the limit set of the associated sample clouds. We first prove a lemma.   
\blm\label{m42}
Let $\WB_n=(U_n,V_n)$ be the componentwise maximum of the sample $N_n=\{\ZB_{n1},\ldots,\ZB_{nn}\}$ from the distribution $\p_n$ on $\rbb^2$. Suppose $n\p_n\{(0,\nf)^2\}\to0$, $n\p_n\{u>0\}\to\nf$, and $n\p_n\{v>0\}\to\nf$.
Then the probability $p_n=\pbb\{\WB\in N_n\}$ that the coordinatewise maximum is a sample point vanishes for $n\to\nf$.  
\elm
\bew Since $\WB_n$ lies in the positive quadrant or $U_n\le0$ or $V_n\le0$, one finds
$$\pbb\{\WB_n\in N_n\}\le\pbb\{N_n((0,\nf)^2)>0\}+\pbb\{N_n(\{u>0\})=0\}+\pbb\{N_n(\{v>0\})=0\}.$$
These binomial probabilities all three vanish for $n\to\nf$. \qed

\begin{figure}[htb]
\centering
\subfloat[]{
\label{sfnd}
\includegraphics[width=0.45\linewidth]{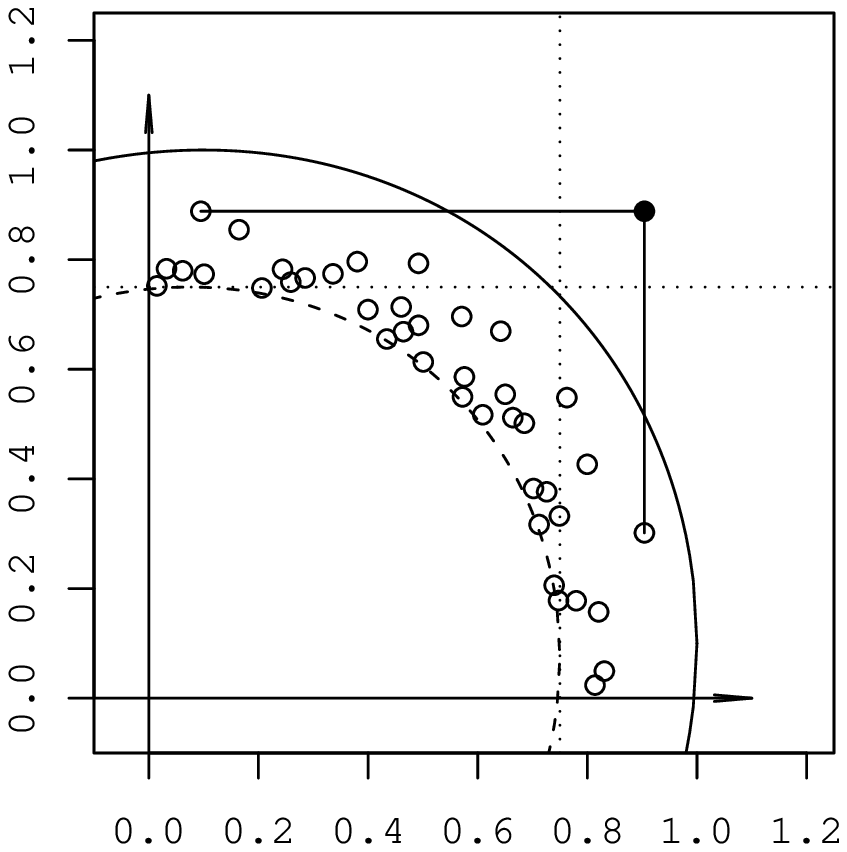}}
\hspace{0.01\linewidth}
\subfloat[]{
\label{sfmetad}
\includegraphics[width=0.45\linewidth]{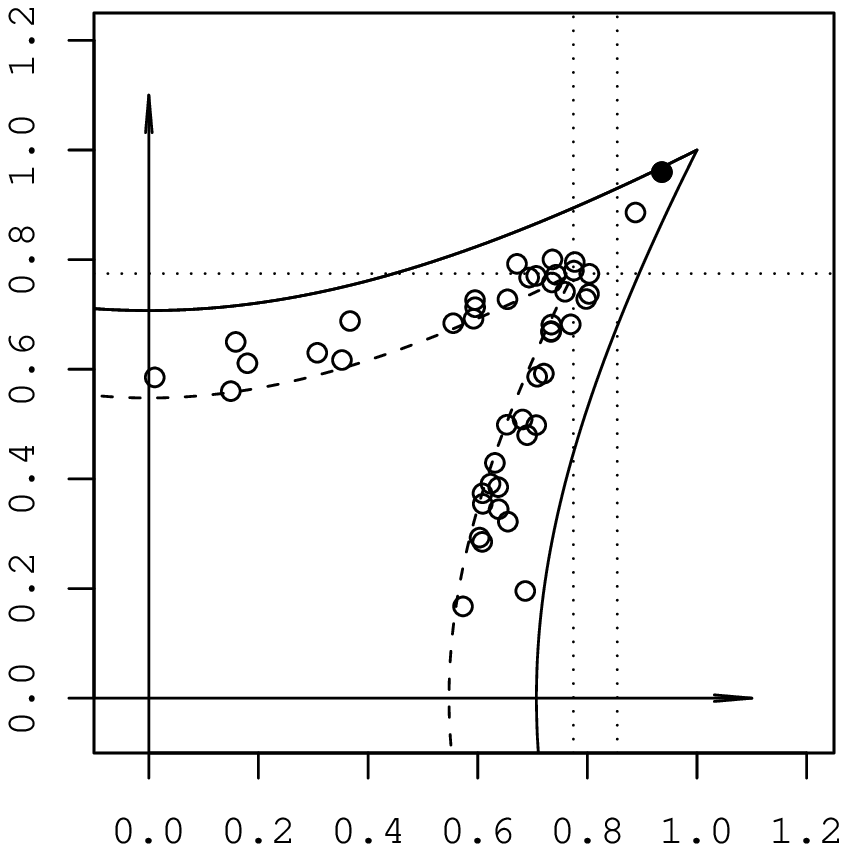}}
\caption{Points at the edge of simulated sample clouds of size $n=10^5$ from (a) a bivariate normal distribution with mean zero and correlation $\r=0.1$ and (b) a bivariate meta-Cauchy distribution with standard normal marginals, in both cases restricted to $(0,\nf)^2$. Sample points are scaled by factor $s_n=\sqrt{2\log n}$. The filled circles indicate the coordinatewise maxima for these samples. The boundaries of the corresponding limit sets $E=\{x^2-2\r xy+y^2\le 1-\r^2\}$ for (a) and $E=\{|x|^2+|y|^2+1\ge3\|(x,y)\|_\nf^2 \}$ for (b) are depicted by solid curves. The sample points inside the dashed curves are not displayed.}
\label{fedge} 
\end{figure}

\bth\label{thmai2}
Let $\XB_1,\XB_2,\ldots$ be independent observations from a continuous df $F$ on $\rbb^d$. Let $D$ be an open bounded star-shaped set which belongs to $\DC_d$. Suppose there exist scaling constants $s_n$ such that the scaled sample clouds $N_n=\{\XB_1/s_n,\ldots\XB_n/s_n\}$ converge onto the closure of $D$. If the bivariate projections $D_{ij}$ of the set $D$ are blunt then $F$ is asymptotically independent.
\eth
\bew For the sake of simplicity assume $d=2$. Let $\bb=(b_1,b_2)$ denote the coordinatewise supremum of $D$. Since $D$ is blunt there exists a $\d>0$ such that the shifted quadrant $(e^{-\d}\bb,\nfb)$ and the set $e^\d D$ are disjoint (see Figure~\ref{sfnd}). Let $n\p_n$ be the mean measure of $N_n$. Then $n\p_n(e^{-\d}\bb,\nfb)\le n\p_n(e^\d D)^c\to0$ and both $n\p_n(D\cap\{v>e^{-\d} b_1\})$ and $n\p_n(D\cap\{u>e^{-\d} b_2\} )$ go to infinity by Definition~\ref{dconv2}. The lemma above applied to the shifted clouds shows that the probability $p_n$ of a record in the sample cloud $N_n$ vanishes for $n\to\nf$. \qed

In order to complete the proof of Theorem~\ref{thmai1}, we now establish a link between the asymptotic shape of the scaled level sets of a light-tailed density and the shape of the limit set of the associated sample clouds. 
\bpr\label{pssc}
Let $\XB_1,\XB_2,\ldots$ denote i.i.d. random vectors from a density $g\in\AC(D)$. Then the sequence of scaled sample clouds $N_n=\{\XB_1/s_n,\ldots,\XB_n/s_n\}$ converges onto the closure of the set $D$ as $n\to\nf$ if the scaling constants $s_n$ are chosen appropriately.
\epr
\bew  The density $g$ is shape equivalent to $f=f_0(n_D)$; see Proposition~\ref{pach}. The function $f_r(\ub)=f(r\ub)/f_0(r)$ also lies in $\HC(D)$, and by rapid variation and monotonicity of $f_0$ for any $M>1$ and $\e>0$ eventually $f_r>M$ on $e^{-\e}D$ and $f_r<1/M$ off $e^\e D$. By the pointwise inequality~\eqref{qACfg} the functions $g_r(\ub)=g(r\ub)/f_0(r)$ satisfy the same inequalities eventually if we replace $\e$ by $2\e$. The measure $\r_r$ with density $g_r$ satisfies $\r_r(D^c)<<\r_r(U)$, $r\to\nf$, for any open set $U$ which intersects $D$ by Lemma~\ref{pAC2}. Choose $s_n$ such that $ns_n^df_0(s_n)=1$ for $n\ge n_0$. Then $\r_{s_n}$ is the mean measure of the sample cloud $N_n$, and from Proposition~\ref{pconv2} the sample clouds $N_n$ converge onto $D$ since their mean measures do.
\qed

\section{Examples}\label{s:examples}

This section illustrates the applicability of Theorems~\ref{thmai3}, \ref{thmai1} and \ref{thmai2} in determining whether a given distribution has asymptotically independent components. We also wish to see whether the conditions are sharp. 

For any open bounded convex set $D$ in the plane whose closure contains the origin, the function $e^{-n_D}$ is integrable and hence $f(\xb)=c_0e^{-n_D(\xb)}$ for suitable $c_0>0$ is a probability density. Along rays it is an exponential function. If $\UB$ is uniformly distributed on $D$ then it is simple to decide when the coordinates $U_1$ and $U_2$ are asymptotically independent, but asymptotic independence of the the coordinates of the vector $\XB$ with density $f $ is a different matter, even in the simple example where $D$ is the intersection of a disk of radius $r=2$ centered in $(1,-1)$ and the open set above the diagonal. The light-tailed examples below are of a more general nature.

\bex{\bf(Rotund-exponential densities)\ }
Let $\XB$ have a continuous homothetic density $f$ with convex shape $D$ and generating function $f_0$. If $D$ has a $C^1$ boundary and $f_0$ varies rapidly, the coordinates are asymptotically independent, and the sample clouds, properly scaled, converge onto the closure of the set $D$ by Theorem~\ref{thmai3} and Proposition~\ref{prd1}. Now assume more: $D$ is rotund, i.e. the boundary $\prl D$ is $C^2$ with positive definite curvature in every point. Also assume that the generating function $f_0$ is asymptotic to a von Mises function $e^{-\j}$. Then $f$ is a so-called rotund-exponential density, see~\cite{Balkema2007}, Sections~9 and~10. If one zooms in onto a boundary point of $D$ so as to distinguish individual sample points, the sample clouds converge to a Gauss-exponential point process: $N_n\imp N$ vaguely on $\rbb^d$. The limit $N$ is a Poisson point process with intensity $g(\ub)=e^{-u_d}e^{-(u_1^2+\cdots+\u_{d-1}^2)/2}/(2\p)^{(d-1)/2}$ if one chooses the normalization appropriately. Weak convergence holds on all halfspaces $\{u_d\ge c_0+c_1u_1+\cdots+c_{d-1}u_{d-1}\}$. The restriction of $g$ to the upper halfspace $\{u_d\ge0\}$ is a probability density. The corresponding vector has independent components. This vector is the limit of the high risk scenarios $\XB^H$, properly normalized, where $\XB^H$ is the vector $\XB$ conditioned to lie in the halfspace $H$, and $H$ moves off to infinity in the sense that $\pbb\{\XB\in H\}\to0$. These results remain valid under certain perturbations (if the density $f=e^{-\j(n_D)}$ is multiplied by a flat function $L$, see Section~11 in \cite{Balkema2007}). Such a perturbation does not affect the asymptotic behaviour of the exponent $\j$, but may affect the limit shape (take $\j=\log^2(1+r)$ with $r^2=x^2+y^2$ and $L=e^\l$ with $\l=\log(1+2r+cx)$, $c\in[-1,1]$. If we choose $c=c(r)=\sin(\log\log r)$ then $L$ still is flat, but the shape of the level sets no longer converges). See Proposition~14.1 in \cite{Balkema2007}. So we see that under the extra conditions on the homothetic density $f$ there are three alternative asymptotic descriptions of large sample clouds from this density. Convergence onto the closure of the set $D$ describes the global behaviour of the sample clouds; weak convergence  in the space $\XC=[-\nf,\nf]^d\sm\{-\nfb\}$ to a Poisson point process whose mean measure is the exponent measure of an extreme value limit law (Gumbel with independent components); weak convergence to the Gauss-exponential Poisson point process $N$ on certain halfspaces describes the local behaviour in boundary points of $D$. \qend
\eex

\bex{\bf (Skew-normal densities)\ }
A symmetric density $g$ satisfies $g(-\xb)=g(\xb)$. It may be transformed into an asymmetric density by multiplication with a positive continuous asymmetric function $\q$ which satisfies $\q(\xb)+\q(-\xb)=2$. The skew-normal distributions $SN(\Omega,\alb)$ introduced in~\cite{Azzalini1996} have a density $f$ which is the product of a centered Gaussian density with covariance $\Omega$ and the function $\xb\mapsto2\Phi(\alb^T\xb)$, where $\Phi$ is the standard normal df and $\alb$ a non-zero linear functional. These densities are log-concave (since $(-\log\Phi)''$ is positive) and hence have convex level sets.

\begin{figure}{htb}
\centering
\subfloat[]{
\label{sfMSN1}
\includegraphics[width=0.45\linewidth]{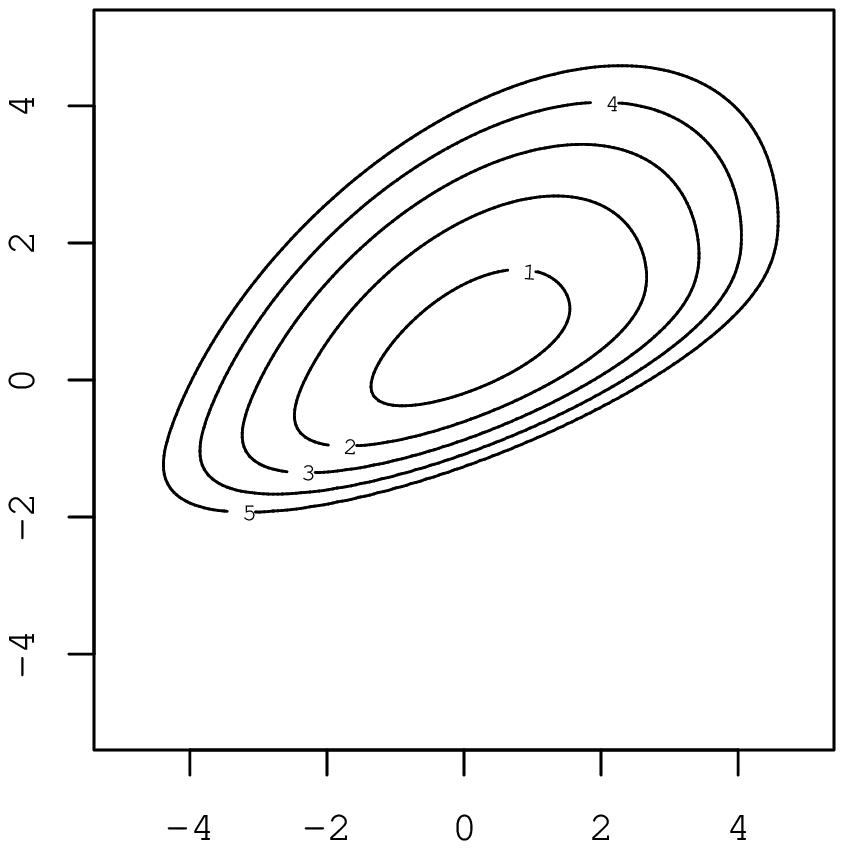}}
\hspace{0.01\linewidth}
\subfloat[]{
\label{sfMSN2}
\includegraphics[width=0.45\linewidth]{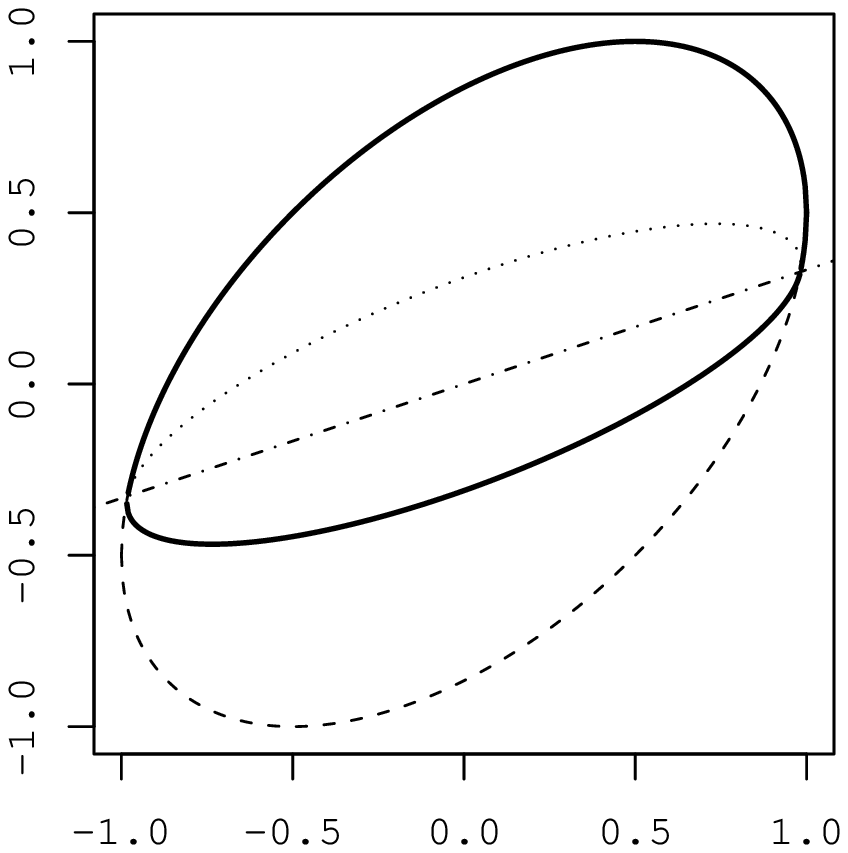}}
\caption{(a) Contours of the density of $SN_2(\alb;\w)$ distribution with $\w=0.5$ and $\alb=(-1,3)$; levels are indicated as powers of $10^{-1}$. (b) Asymptotic shape of the scaled level sets (bold curve); the dash-dotted straight line is given by equation $\alb^T\ub=0$.}
\label{fMSN} 
\end{figure}

We claim that $f\in\AC(D)$ for a convex set $D$ with $C^1$ boundary. Write $f=e^{-h}$. The function $-\log\Phi(t)$ is asymptotic to $t^2/2$ for $t\to-\nf$ and vanishes for $t\to\nf$. Let $\ub_0$ satisfy $\alb^T\ub_0>0$ and $\ub_0^T\Omega\inv\ub_0=1$. Then $h(t\ub_0)/t^2\to1/2$ and
$$\frac{h(t\ub)}{h(t\ub_0)}\to\n(\ub)=\bcs\ub^T\Omega\inv\ub&{\rm for\ }\alb^T\ub\ge0\\\ub^T\Omega\inv\ub+(\alb^T\ub)^2&{\rm for\ }\alb^T\ub<0\ecs\qquad t\to\nf.$$
It follows from Proposition~\ref{prd1} that there is a limit shape: $D=\{\n<1\}$, which is the covariance ellipsoid on the halfspace where $\alb^T\ub$ is positive, and a flattened version of this ellipsoid on the complementary halfspace, see Figure~\ref{fMSN}. The set $D$ is convex. To see that the boundary is $C^1$, choose coordinates such that the underlying Gaussian density is standard, with spherical level sets, and then choose the vertical coordinate in the direction of $\alb$. In these coordinates $D$ agrees with the unit ball $B$ on the upper halfspace, and with the cylinder symmetric ellipsoid $\{x_1^2+\cdots+x_{d-1}^2+(1+c^2)x_d^2<1\}$ for some $c>0$ on the lower halfspace $\{x_d\le0\}$.  For boundary points $\pb$ in the horizontal coordinate plane the tangent plane is vertical: $\pb^T\xb=1$. Theorem~\ref{thmai3} applies. If $\XB$ has a skew-normal distribution then the random variables $\x_1(\XB),\ldots,\x_m(\XB)$ are asymptotically independent whenever the linear functionals $\x_1,\ldots,\x_m$ are linearly independent. Asymptotic independence of the skew-normal distribution has been partially proven in \cite{Lysenko2009} using a direct analytic approach based on Sibuya's condition. \qend
\eex


\bex Densities of the form $f(\xb)=f_0(\|\xb\|_p)$ for $p>1$ have level sets which are balls in $\lp$. Lower dimensional marginals have the same form but with a different generator, since the projection of the $d$-dimensional unit ball on the space spanned by the first $m$ coordinates is the $m$-dimensional unit ball. The  two-dimensional unit ball is blunt for all $p\in[1,\nf)$, and hence vectors $\XB$ with light-tailed densities $f$ as above have asymptotically independent coordinates. However for $p=\nf$, the sup-norm, the unit ball is a cube $C=(-1,1)^d$, and the square is not blunt. The components of $\XB$ are still asymptotically independent, but we need extra work to prove this. The bivariate margins of $f(\xb)=f_0(\|\xb\|_\nf)$ have the same form with a different generator $f_0$ which still is continuous, strictly decreasing and rapidly varying in $\nf$. (The cubic slices are replaced by square slices.) It suffices to consider bivariate densities $f(x,y)=f_0(|x|\lor|y|)$ for continuous, strictly decreasing positive functions $f_0$ on $[0,\nf)$ which vary rapidly in $\nf$. The marginals $f_i$ of $f$, $i=1,2$, are equal by symmetry, and $f_2(y)=2yf_0(y)+2R(y)$ where $R(y)=\int_y^\nf f_0(t)dt<<yf_0(y)$ by rapid variation of $f_0$. Hence $2R(t)/f_1(t)\to0$ for $t\to\nf$, and by l'Hospital's rule also the quotient $\pbb\{X>t,Y>t\}/\pbb\{X>t\}$. Thus Sibuya's condition holds. However, asymptotic independence need not hold if the level sets are only asymptotically cubic, see Example~2 in \cite{Balkema2009a}.
\qend\eex
Another example showing that results for $\HC(D)$ do not need to carry over to $\AC(D)$ is given in the Appendix, Example~\ref{e70}. Our last example illustrates an asymptotically dependent distribution with a density in $\AC(D)$ where bivariate projections of $D$ are non-blunt. 

\bex Let $\ZB$ have a bivariate $t$ density with $\l>0$ degrees of freedom. Transform the marginals to obtain a vector $\XB$ with standard Gaussian components. The new density $g$ is called a meta-$t$ density; see~\cite{McNeil2005} p.193. It has normal marginals but the copula of the elliptic $t$ distribution. The shape of the level sets of the density $g$ converges to the symmetric subset $D=\{u_1^2+u_2^2+\l>(\l+2)\|(u_1,u_2)\|_\nf^2\}$ of the square $(-1,1)^2$, see~\cite{Balkema2009}. Figure~\ref{sfmetad} shows a detail. The set $D$ is not blunt, and the components of $\XB$ are asymptotically dependent since those of $\ZB$ are.\qend \eex

\section{Conclusion}
We have explored conditions for asymptotic independence of the components of a multivariate random vector expressed in terms of the limiting shape of the level sets of the underlying density. A distinction had to be made between light and heavy tails. For light-tailed densities, the limiting shape of level sets is essential in determining whether asymptotic independence holds. In contrast, for heavy-tailed densities the (limiting) shape of level sets is irrelevant as long as the shape intersects the positive orthant. In the light-tailed case there is a simple sufficient condition for asymptotic independence of two components in terms of the corresponding bivariate projection of the shape. This subset of the plane has to be blunt. The more delicate question of the relation between shape and asymptotic independence when the bivariate projection is not blunt will be treated in a future publication. 

Asymptotic dependence is a basic concern in multivariate risk analysis. The light-tailed densities studied in this paper have the property that sample clouds will have the same shape as the level sets of the density asymptotically. For sample clouds persistence of the shape, as the number of sample points increases, opens the possibility of using the shape to construct densities over the whole space. This makes it possible to estimate probabilities of regions far out which contain only a few or no sample points.

\subsubsection*{Acknowledgments} We are thankful to Paul Embrechts for drawing our attention to the problem considered in this paper and for several useful discussions.
\appendix
\numberwithin{equation}{section}
\numberwithin{theorem}{section}
\numberwithin{example}{section}
\numberwithin{subsection}{section}

\section{Appendix}

\subsection{Supplementary results}

\bpr\label{prec}
Let $\XB_1,\XB_2,\ldots$ be independent observations from the continuous df $F$ on $\rbb^2$. The probability of a record amongst the first $n$ observations goes to zero if and only if $F$ is asymptotically independent.\epr
\bew We may assume that $F$ is a copula. Set $c_n=nP(1-1/n,1-1/n)$. Then the Poisson approximation gives a probability $p_n=c_ne^{-c_n}e^{-(1-c_n)}e^{-(1-c_n)}=c_ne^{c_n-2}$ to the event: among the first $n$ observations there is exactly one in the complement of $[0,1-1/n]^2$, and that observation lies in the square $(1-1/n,1]^2$.. In case of asymptotic dependence $c_{k_n}\to c>0$ for some subsequence, and hence the probability of a record in a sample of size $k_n$ will exceed $ce^{c-2}/2$ eventually. Conversely, asymptotic independence implies $nP(1-M/n,1-M/n)\to0$ for any $M>1$, whereas the marginals satisfy $nP(0,1-M/n)=M$. As in Lemma~\ref{m42} the probability of a record vanishes.\qed

\bpr\label{pau1}
Suppose $f\in\HC(D)$ and the generator $f_0$ varies regularly or satisfies $O$-variation 
\bql{q3orv}r_{n+1}\sim r_n\to\nf\imp f_0(r_{n+1})\sim f_0(r_n).\eql
Then $g$ is shape equivalent to $f$ if and only if $g\sim f$.
\epr
\bew Condition~(\ref{q3orv}) implies that for every $\e>0$ there exists $\d>0$ such that $f_0(e^\e r)>e^{-\d}f_0(r)$ for all $r>0$. Let $n_D(\xb_n)\to\nf$. Suppose $g(\xb_n)=c_n=f_0(r_n)$. Then for any $\e>0$ the point $\xb_n$ eventually lies in the ring $e^\d r_nD\sm e^{-\d}r_nD$ on which $f$ fluctuates by a factor at most $e^\e$. To show the converse, suppose $f_0$ satisfies \eqref{q3orv} and let $f(\xb_n)=f_0(r_n)$ so that $g(\xb_n)\sim f(\xb_n)=f_0(r_n)\sim f_0(r_{n+1})$ for $r_{n+1}\sim r_n\to\nf$. Then for any $\e,\e_1>0$ eventually $$(1-\e_1)f_0(e^{\e}r_n)\le (1-\e_1)f_0(r_{n+1})\le g(\xb_n)\le (1+\e_1)f_0(r_{n+1})\le(1+\e_1)f_0(e^{-\e}r_n),$$
and since $\e_1$ is arbitrary, we have $f(e^\e\xb_n)\le g(\xb_n)\le f(e^{-\e}\xb_n)$ for $n\ge n_0$ as required for shape equivalence by \eqref{qag1}.\qed

\subsection{Counterexamples}

This section contains counterexamples mentioned in the main text.

\bex\label{e60}
A density $f\in\HC(D)$ may have spherical level sets without exhibiting spherical symmetry. This will be the case if $D$ is an off-center ball with $n_D(\xb)=\sqrt{\|\xb\|_2^2+\betab^T\xb}-\betab^T\xb$ for some $\betab\in\rbb^d$. To be star-shaped the origin has to lie in $D$, or be a boundary point. In the latter case the set $\{f>0\}$ is a ball (if $\{f_0>0\}$ is a bounded interval) or an open halfspace. \qend
\eex

\bex\label{e30}
For any $\e>0$ there exists a bounded open star-shaped set $D$ which contains the origin, whose closure is the cube $K=[-1,1]^d$, and such that the volume of $D$ is small, $|D|<\e$. To construct such a star-shaped set take a dense sequence $\xb_n$ on the boundary $\prl K$, and define $U\ss\prl K$ as the union of open disks with center $\xb_n$ and radius $\e_n$, where $\e_n\to0$ so fast that the area of $U$ is $\e/2$. Now let $D$ be the union of an open centered ball with volume $\e/2$ and the set of all points $r\ub$ with $0<r<1$ and $\ub\in U$. \qend
\eex

\bex \label{e40} Consider a continuous strictly positive density $f$ on $\rbb^3$ whose level sets $\{f>c\}=r_cD$ all have the same shape. The set $D$ is convex, even rotund, and the function $c\mapsto r_c$ is continuous and strictly increasing. The marginal densities are not necessarily all unimodal.

Let $f_0$ be the uniform density on the tetrahedron $T_0$ with two vertices in the horizontal plane, say $\eb_1$ and $\eb_2$. The other two vertices are $\eb_3/m$ and $-\eb_3$. Here $m$ is a positive integer to be chosen later. The marginal along the vertical axis along the base vector $\eb_3$ has a continuous density $g_0$ on $(-1,1/m)$ which is parabolic on the interval $(-1,0)$ and on $(0,1/m)$, vanishing in the endpoints of the interval and with a maximum in the origin.
The shifted tetrahedron $T=T_0+\eb_3/2-(\eb_1+\eb_2)/8$ contains the origin as interior point. Its vertical marginal density $g$ is $g_0$ shifted upwards over $1/2$ and has its maximum in $1/2$. The vertical marginal density $\tilde g$ of the uniform distribution on the tetrahedron $T/2$ has its maximum in $1/4$. The fair mixture of the uniform distribution on $T$ and $T/2$ has a density $\hat f$ whose marginal $\hat g=(g+\tilde g)/2$ is not convex if $m$ is large (since the left derivative of $\tilde g$ in $1/4$ is large). Now choose rotund sets $D_n$ converging to $T$ and densities $f_n$ converging to $\hat f$ which satisfy the conditions of the example. If infinitely many of the vertical marginals $g_n$ were unimodal then the limit $\hat g$ would be. We conclude that eventually $g_n$ is not unimodal.\qend
\eex

\bex\label{eAAC} In view of Lemma~\ref{prve} we see that $qg$ is shape equivalent to $g$ for a continuous function $q=e^\j$ precisely if $|\j(\xb)|\le\c(\|\xb\|)$ for a function $\c(r)<<r^2$. Here are some examples of functions $h=qg$ which are shape equivalent to the standard normal density $g$ on the plane, for a continuous positive function $q$. One may take $q$ to be one of the following functions
$1+|x|$, $(1+r^2)^m$ with $m\ge1$, $e^r$, $e^{x-|y|^{3/2}}$ where $r^2=x^2+y^2$. These functions may be multiplied with a function like $\exp(\sin\p e^{x^2}\sin\p e^{y^6})$ which fluctuates rapidly but is weakly asymptotic to a constant. The level sets $\{h>c\}$ then will look locally like a shore with many small islands, and lakes, even though the sets are asymptotic to disks $\{x^2+y^2<r^2\}$ with $r=\sqrt{2\log(1/c)}$.
\qend
\eex

\bex\label{e70}
Let $D$ be the open triangle with vertices $(1,1)$, $(-1,0)$ and $(0,-1)$. It contains the origin. Let $f\in\AC(D)$ have convex level sets and be shape equivalent to $g=e^{-n_D}$. It is possible that $f$ is asymptotically independent. Suppose the function $g$ has triangular level sets $\{g>e^{-t}\}=tD$. Let $f$ have level sets $\{f>e^{-t}\}=D_t=tD\sm\{x+y\ge 2t-\sqrt t\}$. Then $D_t/t$ agrees with the triangle $D$ except that the extreme top has been sliced off. If we choose $\pb(t)=(t,t)-(\sqrt t,\sqrt t)/2$, then 
$$h_t(\ub)=f(\pb(t)+\ub)/f(\pb(t))\to h(\ub)\qquad t\to\nf,$$
where $\{h>e^{-t}\}=C+(t,t)$ for the halfspace $C=\{u+v<0\}$. Let $\r_t$ have density $h_t$ and let $\r$ have density $h$. Then $\r_t\to\r$ weakly on $[0,\nf)^2$. Since $\r[0,\nf)^2$ is finite and $\r$ gives infinite weight to the halfspaces $\{y\ge0\}$ and $\{x\ge0\}$, Sibuya's condition holds by Proposition~\ref{p2sa} with curve $\pb(t)$, $t\ge1$.  \qend
\eex


\bibliography{bibliographyAI}{}
\bibliographystyle{plain} 


\end{document}